\documentclass[12pt]{article}
\usepackage[ansinew]{inputenc}
\usepackage{pst-all,epsfig}
\usepackage{amsmath}
\usepackage{amsfonts}
\usepackage{float}
\usepackage{graphicx, graphics}
\usepackage{amssymb}
\usepackage[left=2.5cm,top=2cm,right=2.5cm,bottom=2.5cm]{geometry}
\newcommand{\fin}{\hfill $\Box$}

\usepackage[USenglish]{babel}

\usepackage{graphics}
\newcommand{\Rn}{\mathbb{R}^{n}}

\newcommand{\inte}{\operatorname*{int}}

\newtheorem{lemma}{Lemma}
\newtheorem{theorem}{Theorem}

\newtheorem{remark}{Remark}

\newtheorem{conjecture}{Conjecture}

\title{Sections and projections of nested convex bodies}
\author{
I. Gonz\'alez-Garc\'ia$^1$, J. Jer\'onimo-Castro$^{2}$ \\ E. Morales-Amaya$^{3}$, and D. J. Verdusco-Hern\'andez$^{4}$ \\ 
\small{$^{1,2}$Facultad de Ingenier\'ia}\\
\small{Universidad Aut\'onoma de Quer\'etaro, M\'exico}\\
\small{$^{3,4}$Facultad de Matem\'aticas-Acapulco,}\\
\small{Universidad Aut\'onoma de Guerrero, M\'exico}\\
}

\date{\small{\today}}

\begin{document}
 \maketitle
\begin{abstract}  
One of the most important problems in Geometric Tomography is to establish properties of a given convex body if we know some properties over its sections or its projections. There are many interesting and deep results that provide characterizations of the sphere and the ellipsoid in terms of the properties of its sections or projections. Another kind of characterizations of the ellipsoid is when we consider properties of the support cones. However, in almost all the known characterizations, we have only a convex body and the sections, projections, or support cones, are considered for this given body. In this article we proved some results that characterizes the Euclidean ball or the ellipsoid when the sections or projections are taken for a pair of nested convex bodies, i.e., two convex bodies $K$, $L$ such that $L\subset\text{int}\, K.$ We impose some relations between the corresponding sections or projections and some apparently new characterizations of the ball or the ellipsoid appear. We also deal with properties of support cones or point source shadow boundaries when the apexes are taken in the boundary of $K$.
\end{abstract}
 
\section{Introduction}
One of the most important problems in Geometric Tomography is to establish properties of a given convex body if we know some properties over its sections or its projections. A very interesting topic was initiated by C. A. Rogers \cite{Rogers} who proved: 

\emph{Let $K\subset\mathbb R^n$ be a convex body, $n\geq 3$, and $p$ be a point in $\mathbb R^n$. If all the $2$-dimensional sections, i.e., the intersections of $K$ with planes of dimension 2, through $p$ have a centre of symmetry, then $K$ has also a centre of symmetry.} 

Rogers also conjectured that if $p$ is not the centre of $K$ then it is an ellipsoid. This conjecture was first proved by P. W. Aitchison, C. M. Petty, and C. A. Rogers \cite{false_centre} and is known as the False Centre Theorem. A simpler proof was given by L. Montejano and E. Morales-Amaya  in \cite{Falso_MM}. The following is a remarkable theorem due to S. Olovjanishnikov \cite{Olovja}: 

\emph{A convex body $K\subset\mathbb R^n$, $n \geq 3$, is a solid ellipsoid provided all hyperplanar sections of $K$ that divide the volume of $K$ in a given ratio $\mu\neq 1$ have a centre of symmetry.} 

We say that two sets $A$ and $B$ in $\mathbb R^n$ are affinely equivalent if there is an invertible affine transformation
$f:\mathbb R^n\longrightarrow \mathbb R^n$ such that $A=f(B).$ A result similar to the False Centre Theorem was proved by L. Montejano \cite{Montejano1} which combines the method of P. Mani \cite{Mani} and some results of algebraic topology: 

\emph{Let $K\subset\mathbb R^n$ be a convex body, $n\geq 3$, and $p\in\emph{int}\, K$. If all sections of $K$ by hyperplanes through $p$ are affinely equivalent, then $K$ is either a solid ellipsoid or has centre of symmetry at $p$.} 

Recently, L. Montejano proved a similar result \cite{Montejano2} where he consider affinely equivalent projections instead of affinely equivalent sections. 

A different kind of characterizations of the ellipsoid is when we consider properties of the support cones. Let $K\in\mathbb R^n$, $n\geq 2$, be a convex body, i.e., a compact and convex set with non-empty interior. Given a point $x\in\mathbb R^n\setminus K$ we denote the cone generated by $K$ with apex $x$ by $\text{C}(K,x)$, i.e., $\text{C}(K,x):=\{x+\lambda (y-x): y\in K,\, \lambda\geq 0\}.$ The boundary of $\text{C}(K,x)$ is denoted by $\text{S}(K,x)$, in other words, $\text{S}(K,x)$ is the support cone of $K$ from the point $x$. We denote the \emph{graze} of $K$ from $x$ by $\Sigma(K,x)$, i.e., $\Sigma(K,x):=\text{S}(K,x)\cap \partial K.$ In \cite{Marchaud} A. Marchaud proved the following: 

\emph{Let $K\subset \mathbb R^3$ be a convex body and $H$ be a plane which is either disjoint from $K$ or meets $K$ at a single point. Then $K$ is a solid ellipsoid if for every point $x\in H\setminus K$, the graze $\Sigma(K,x)$ contains a planar convex curve $\gamma$ such that $\emph{conv}\, \gamma \cap \emph{int}\, K \neq\emptyset$}. 

This result is also true if the apexes of the cones are points at infinity, i.e., the grazes, in this case called \emph{shadow boundaries} are obtained by intersections of $\partial K$ with circumscribed cylinders. In this case is more convenient to think that the body is illuminated by light rays parallel to a given direction $u\in\mathbb S^{n-1}$, and the shadow boundary is precisely the border between the light part and the dark part of $\partial K$. The shadow boundary of $K$ with respect the direction of illumination $u$ is denoted by $\Lambda (K,u)$. The firs proof that a convex body is an ellipsoid if and only if every shadow boundary lies in a hyperplane is due to H. Blaschke \cite{Blaschke}. However, if we change the plane $H$ by a surface $\Gamma$ which encloses to $K$ it is not known whether if $K$ is an ellipsoid or not, even in the case where $\Gamma$ is the boundary of a convex body that contains $K$ in its interior. We suspect the following is true.

\begin{conjecture}\label{penumbras_planas} Let $L\subset \emph{int}\, K\subset\mathbb R^n$ be convex bodies such that for every point $x\in\partial K$ it holds that $\Sigma(L,x)$ lies in a hyperplane. Then $L$ is an ellipsoid.
\end{conjecture}

Another characterization of ellipsoids with respect to the properties of support cones is the following due to P. Gruber and T. \'Odor \cite{Gruber_Odor}: 

\emph{Let $K\subset \mathbb R^n$, $n\geq 3$, be a convex body whose boundary is a hypersurface of class $C^4$ with positive Gauss curvature. If for every exterior point $x$ of $K$ which is sufficiently close to $K$, the cone $\emph{C}(K,x)$ has an axis of symmetry, then $K$ is an ellipsoid.}

In this article we proved some results that characterizes the Euclidean ball or the ellipsoid when the sections or projections are taken for a pair of nested convex bodies, i.e., two convex bodies $K$, $L$ such that $L\subset\text{int}\, K.$ We impose some relations between the corresponding sections or projections and some apparently new characterizations of the ball or the ellipsoid appear. We also deal with properties of support cones or point source shadow boundaries when the apexes are taken in the boundary of $K$.

\section{By properties of sections or orthogonal projections}
We give first the definition of a very important convex
body. For the given convex body $K$, the floating body of $K$ is defined as follows: for a given
positive number $\delta < \text{vol}_n(K)$, the \emph{floating
body} denoted as $K_\delta$ is defined as the intersection of all
the closed half-spaces which cut off from $K$ a cap with volume
$\delta$. For instance, if $K$ is a Euclidean disc with unit area
in the plane and $\delta <1/2$, then $K_{\delta}$ is a disc
concentric with $K$. Notice that the midpoint of every chord of
$K$ and tangent to $K_{\delta}$ belongs to $K_{\delta}.$ This
property holds for every convex body and its floating body
$K_{\delta}$ (see for instance \cite{SW1}).

\begin{figure}[H]
    \centering
    \includegraphics[width=.43\textwidth]{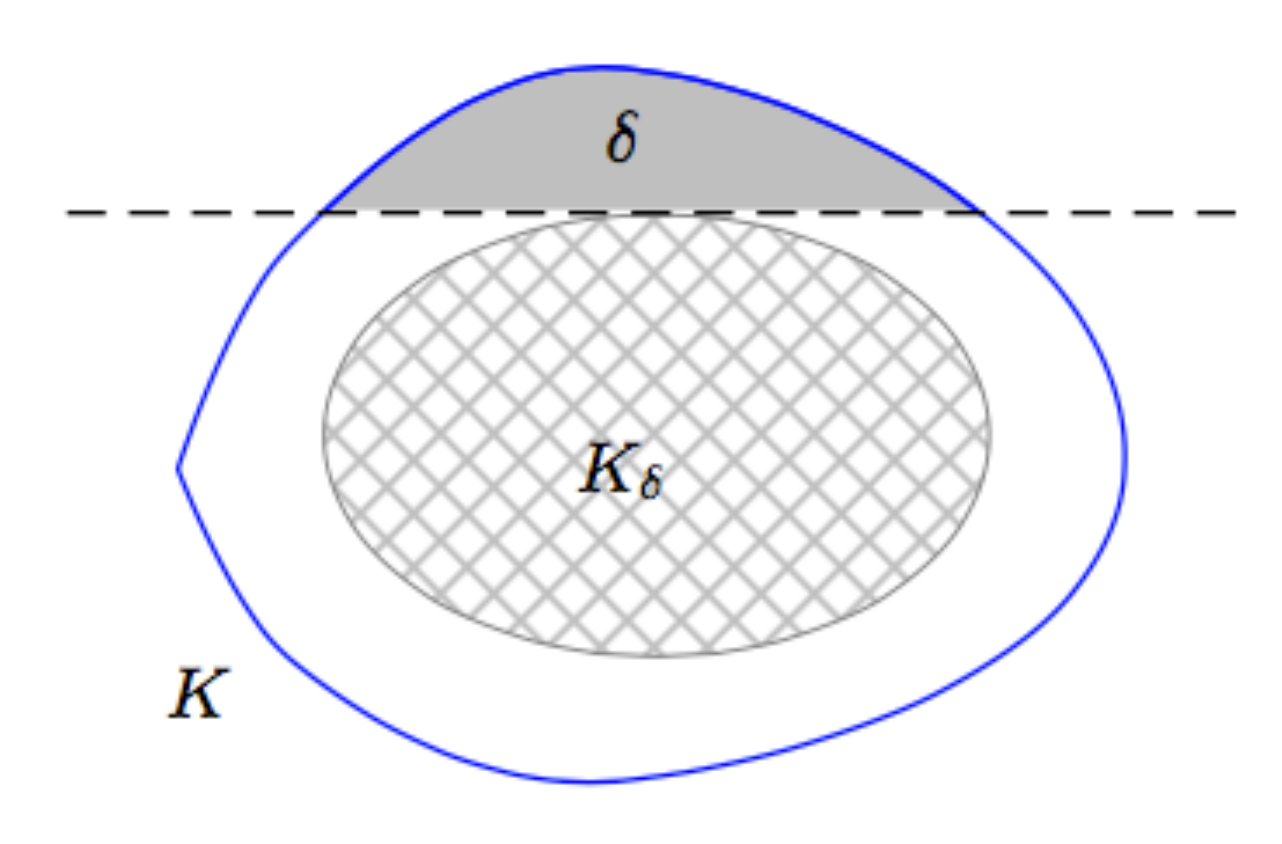}
    \caption{The floating body}
    \label{cuerpos}
\end{figure}

The following result is related to the notions of floating body and a special set of hyperplanes called an $n$-cycle of hyperplanes (see for instance \cite{JMM}): Let $\delta :\mathbb{S}^n\longrightarrow\mathbb{R}$ be a continuous function such that $\delta (-x)=-\delta (x)$. We denote by $\mathcal P$ the following set of hyperplanes in $\mathbb{R}^{n+1}$:
$$\mathcal P =\{ H_y = \{x\in\mathbb{R}^{n+1}|\langle x,y\rangle = \delta (y)\}\}_{y\in\mathbb{S}^n},$$ and we say that $\mathcal P$ an $n$-cycle of hyperplanes. If for a given convex body $K$, it holds that $H_p\cap H_q\cap \text{int}\, K\neq\emptyset$ for every two hyperplanes $H_p,H_q\in\mathcal P$, we say that $\mathcal P$ is an $n$-cycle of planes for $K$. In \cite{Larman}, D. Larman, L. Montejano, and E. Morales proved some results which characterizes the ellipsoid (in $\mathbb R^3$) in terms of the properties of the sections, cut by the members of a $2$-cycle of planes, in a given convex body. 

We have the following interesting lemma.

\begin{lemma}\label{nciclo}
Let $K\subset\mathbb{R}^{n+1}$ a convex body and let $\ell\subset K$ be any $(n-1)$-dimensional section of $K$. If $\mathcal P$ is an $n$-cycle of hyperplanes for $K$ then there exists $\Pi\in\mathcal P$ such that $\ell\subset\Pi$. 
\end{lemma}

\emph{Proof.} For every unitary vector $u\in\mathbb{S}^n$, we may consider that $\delta(u)$ is the distance with sign from $H_u$ to the origin $O$. Recall that $H_u$ is the plane in $\mathcal P$ that is orthogonal to $u$.  Let $\mathcal P_{\ell}$ be the subset of $\mathcal P$ such that for every $H_y\in\mathcal P_{\ell}$ it holds that $y$ is orthogonal to $\text{aff}\, \ell$. Denote by $S(\ell)=\{\upsilon\in\mathbb S^n: \upsilon \bot \text{aff}\, \ell\}.$ Let $\delta_{\ell}$ be the restriction of $\delta$ to the set $S(\ell)$. Let $O'$ be the centroid of $\ell$, and let $\sigma: S(\ell)\longrightarrow\mathbb R$ be the function such that $\sigma (x)$ is the distance with sign from $O'$ to $H_x$. Since $\delta_{\ell}$ is a continuous function we have that $\sigma$ is also a continuous function. Moreover, $\sigma(-x)=-\sigma(x)$, i.e., $\sigma$ is an odd function. By the Borsuk-Ulam theorem (see for instance \cite{Matousek}), there exists $x_0\in S(\ell)$ such that $\sigma(x_0)=0$, this means that the plane $H_{x_0}$ contains to $\ell$. \fin

The following is a more convenient reformulation, four our purposes, of the theorem of Olovjanishnikov \cite{Olovja}.

\textbf{Theorem O.} Let $K$ and $L$ be convex bodies in $\mathbb R^n$, with $L\subset\text{int}\, K$, such that every hyper-section of $K$ tangent to $L$ is centrally symmetric and its centre belongs to $\partial L$, then $K$ and $L$ are homothetic and concentric ellipsoids. 

\begin{theorem}\label{secciones_flotacion}
Let $K, L\subset \mathbb R^n$, $n\geq 3$, be two convex bodies with $L\subset \emph{int}\, K$. Suppose that one of the following properties is satisfied:
\begin{enumerate}
\item [(a)] There is an $(n-1$)-cycle of planes $\mathcal P$ for $L$ such that for every $H\in\mathcal P$ it holds that $H\cap L$ is a floating body of $H\cap K$. 
\item [(b)] There is a sphere $\mathcal S$ contained in $\emph{int}\, L$, such that for every supporting hyperplane $H$ of $\mathcal S$ we have that $H\cap L$ is a floating body of $H\cap K$.
\end{enumerate}
Then $K$ and $L$ are homothetic and concentric ellipsoids.
\end{theorem}

\emph{Proof.} First we suppose that (a) holds. Let $x$ be any point in $\partial L$ and let $\ell$ be any supporting $(n-2)$-dimensional plane through $x$. By Lemma \ref{nciclo} we have that there exists a hyperplane $H\in\mathcal P$ such that $H\cap L$ is a floating body of $H\cap K$. Since $\ell\subset H$ is a supporting $(n-2)$-dimensional plane of $H\cap L$ through $x$, we have that $x$ is the centroid of $\ell\cap K$ (see for instance \cite{SW1}). Because $\ell$ is any $(n-2)$-dimensional supporting plane of $L$ through $x$, we have that $x$ is the centre of symmetry of the section of $K$ tangent to $L$ at $x$ (see Theorem 5.6.19 in \cite{Gro}). By Theorem O we have that $K$ and $L$ are homothetic  and concentric ellipsoids. 

Now, suppose that (b) holds. Let $x$ be any point in $\partial L$ and let $\ell$ be any supporting $(n-2)$-dimensional plane through $x$. Consider one of the two supporting hyperplanes of $\mathcal S$ that contain to $\ell$, to say $H$. From this point we proceed exactly as in the proof of (a). Hence $K$ and $L$ are homothetic  and concentric ellipsoids. \fin

With respect to projections we have the following conjecture.

\begin{conjecture}\label{proyecciones_flotacion}
Let $K, L\subset \mathbb R^n$, $n\geq 3$, be two convex bodies with $L\subset \emph{int}\, K$ a strictly convex body (i.e., with no segments on its boundary). Suppose for every $u\in\mathbb S^{n-1}$ the orthogonal projection of $L$, in direction $u$, is a floating body of the orthogonal projection of $K$. Then $K$ and $L$ are homothetic and concentric ellipsoids.
\end{conjecture}

As evidence of the veracity of this conjecture we have a proof for the case $n=3$.

\emph{Proof for $n=3$.} Let $x\in\partial L$ be any point of $L$, i.e., there exists a support plane $H$ of $L$ such that $H\cap L:=\{x\}.$ Let $u\in\mathbb S^{2}$ be any vector parallel to $H$ and denote by $\pi_u$ the orthogonal projection onto the subspace $u^{\bot}$. The projection of $\pi_u(H)$ is a supporting line of $\pi_u(L)$ through the point $x':=\pi_u(x)$ (see Fig. \ref{proyecciones_flotacion}). 

\begin{figure}[H]
    \centering
    \includegraphics[width=.92\textwidth]{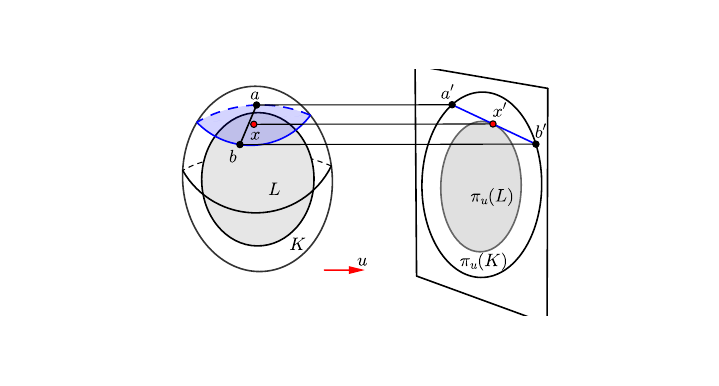}
    \caption{The point $x$ projects onto the midpoint of $[a',b']$}
    \label{proyeccion_flotacion}
\end{figure}

Since $\pi_u(L)$ is the floating body of $\pi_u(K)$ we have that $x'$ is the midpoint of the segment $[a',b']:=H\cap \pi_u(K)$ (see for instance \cite{SW1}). This is true for any vector $u$ parallel to $H$, we then have that $x$ is the centre of symmetry for $H\cap K$. Now, since this last is true for every point $x\in\partial L$, we apply Theorem O and conclude that $K$ and $L$ are homothetic and concentric ellipsoids.  \fin

The proof of the general case of this conjecture depends upon the following conjecture.

\begin{conjecture}\label{centroides}
Let $K\mathbb\subset \mathbb R^n$, $n\geq 3$, be a convex body and let $q$ be a point in its interior. If for every $u\in\mathbb S^{n-1}$ we have that $\pi_u(q)$ is the centroid of $\pi_u(K)$, then $q$ is the centre of symmetry of $K$.
\end{conjecture}

We were not able to prove it or to find this in the literature, neither we have a counterexample.

\begin{remark}
Let $K$ and $L$ be convex bodies in $\mathbb R^3$, with $L\subset\emph{int}\, K$. If for every hyperplane $H$, the projection of $L$ onto $H$ is the floating body of the projection of $K$, then $K$ and $L$ are homothetic and concentric ellipsoids and then $L$ is a floating body of $K$. However, the converse is not true, i.e., if $L$ is the floating body of $K$, then it is not necessarily true that every projection of $L$ is the floating body of the projection of $K$.
\end{remark}

\begin{remark}
Let $K\subset\mathbb R^n$ be a convex body with boundary of class $C^1$. A billiard ball is a point in $K$ which moves along a straight line in the interior of $K$ until it hits $\partial K$ where it is reflected in the usual way: the angle of reflection, with respect to the normal line at the point, is equal to the angle of incidence and the incidence line, the normal line, and the reflection line, are coplanar. A convex body $L$ in the interior of $K$ is a convex caustic  if any trajectory which touches $L$ once, touches it again after each reflection.
In \cite{Gruber}, P. Gruber proved the following: \emph{Let $K\subset\mathbb R^n$ be a convex body with boundary of class $C^1$, let $L$ be a convex body in $\text{int}\, K$, and let $r\in\{2,3, \ldots ,n-1\}$. Then $L$ is a caustic of $K$ if and only if for every $r$-dimensional plane $H$, the orthogonal projection of $L$ onto $H$ is a caustic of the orthogonal projection of $K$.} In some sense, Conjecture \ref{proyecciones_flotacion} is similar to Gruber's theorem. Indeed, in \cite{Gruber2} Gruber proved that only ellipsoids have convex caustics.
\end{remark}

\section{By properties of the support cones}
Let $K\in\mathbb R^n$, $n\geq 2$, be a convex body, i.e., a compact and convex set with non-empty interior. We say that a hyperplane $\ell$ is a support hyperplane of $K$ if \ $\ell\cap K\neq\emptyset$ \ and $K$ is contained in one of the half-spaces bounded by $\ell$. In the Euclidean plane we have the following: for every $t\in[0,2\pi]$ denote by $\ell(t)$ the support line of $K$ with outward normal vector $u(t)=(\cos t,\sin t),$ and let $p(t)$
denote the distance with the sign from the origin $O$ to
$\ell(t).$ A well known result in the geometry of convex sets is that a body is determined in a unique way if the support function is given (see for instance \cite{Val}). Particularly, a convex body is a disc centred at the origin if its support function is constant.

In \cite{Jero_McAlly}, J. Jer\'onimo-Castro and T. McAllister proved the following: \emph{Let $\mathcal C\subset\mathbb R^n$, $n\geq 3$, be a convex solid cone. If for every $y\in\mathbb R^n$ the intersection of $(\partial \mathcal C)\cap (y-\partial \mathcal C)$, when non empty, lies in a hyperplane, then $K$ is an elliptic cone, i.e., every $(n-1)$-dimensional bounded section of $K$ is an $(n-1)$-dimensional ellipsoid.} The following theorem is similar to this mentioned theorem, however, there is more order in the position of the apexes of the cones.

\begin{theorem}\label{secciones_coinciden}
Let $K$ be a convex body contained in the interior of the unit ball of $\mathbb R^n$, $n\geq 2$, such that for every unit vector $x\in\mathbb S^{n-1}$ we have that $S(K,x)\cap S(K,-x)$ is contained in the subspace $x^{\bot}.$ Then $K$ is a Euclidean ball.
\end{theorem}

\emph{Proof.} We first note that the origin $O$ is contained in the interior of $K$. Suppose this is not the case, then there exists a hyperplane $H$ which separates $O$ from $K$ (it could happen that $O\in H$). Let $z\in \mathbb S^{n-1}$ be the vector orthogonal to $H$ which is contained in the opposite half-space bounded by $H$ that $K$. Clearly, the intersection $S(K,z)\cap S(K,-z)$ cannot be contained in $z^{\bot}$. This contradiction shows that $O$ must be in the interior of $K$.

Now we will prove that the central projections of $K$ onto $x^{\bot}$, from the points $x$ and $-x$, coincide for every $x\in\mathbb S^{n-1}.$ In order to do this, let $L$ be any $2$-dimensional plane through the points $x$ and $-x$, and consider the section $K_L:=K\cap L$. Consider the support lines of $K_L$, from $x$ and $-x$, that touch $K_L$ in the same half-plane bounded by the line through $x$ and $-x$. The intersection point $q$ of these support lines must be, by hypothesis, in the hyperplane $x^{\bot}$. This argument proves the assertion.

Now, we continue the proof by induction.
\begin{enumerate}
\item For $n=2$. Let $x_0\in\mathbb S^1$ be such that the angle $\beta$ between one of the support lines of $K$ from $x_0$ and the line generated by $x_0$ is an irrational multiple of $\pi$ (as shown in Fig. \ref{sucesion_circulo}). The existence of such $x_0$ is guaranteed by the following argument: if for every $x\in\mathbb S^1$ the associated angle $\beta$ between the corresponding support line and the vector $x$ is equal to a constant value $\beta_0$, then the support function of $K$ is also constant and consequently, $K$ is a disc (see \cite{Val}). Otherwise, there is an interval $[\beta_1,\beta_2]$ such that $\beta$ takes every value in this interval in a continuous way. Obviously, there is a value in this interval such that $\beta$ is an irrational multiple of $\pi$. Now, let $x_1$ be the intersection between this support line from $x_0$ and $\mathbb S^1$. Denote by $K_{x_0}:=C(K,x_0)\cap x_0^{\bot}$, i.e., $K_{x_0}$ is the central projection of $K$ onto $x_0^{\bot}$ from the point $x_0$. By hypothesis $K_{-x_0}=K_{x_0}$, hence the angle between the corresponding support line of $K$ from $-x_0$ and the line generated by $-x_0$ is also $\beta$. 

\begin{figure}[H]
    \centering
    \includegraphics[width=0.51\textwidth]{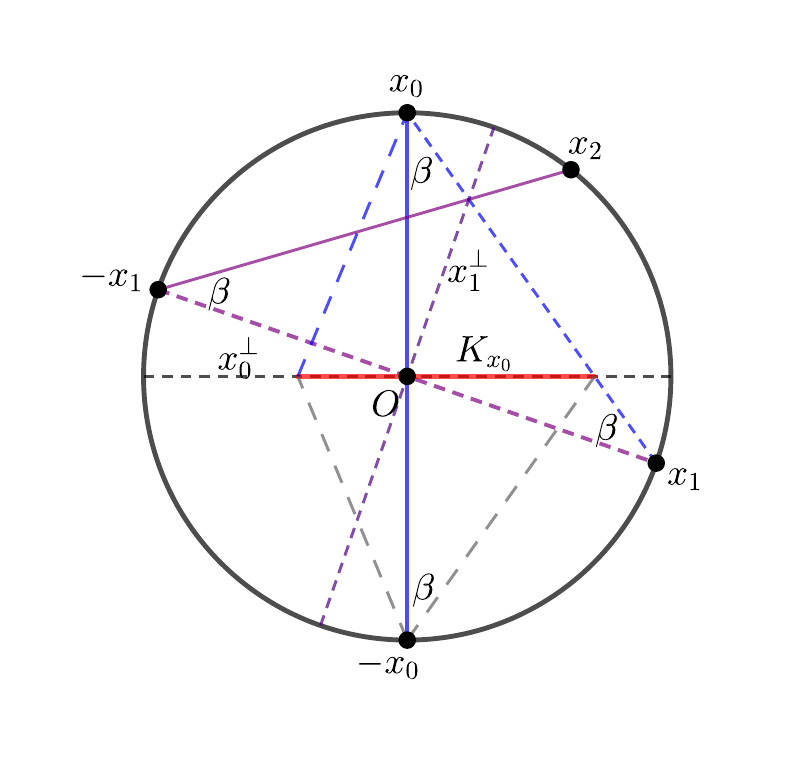}
    \caption{}
    \label{sucesion_circulo}
\end{figure}

Since triangle $\triangle x_0Ox_1$ is isosceles, we have that $\measuredangle Ox_1x_0=\beta$. Let $z$ be the intersection between the segment $[x_0,x_1]$ and the line $x_1^{\bot}.$ By the hypothesis of the theorem we have that the angle between the segments $[x_1,-x_1]$ and $[-x_1,z]$ is equal to $\beta$. Now, let $x_2$ be the point in $\mathbb S^1$ (different of $-x_1$) that is collinear with $-x_1$ and $z$. Since $\measuredangle x_2Ox_1$ is a central angle which intersects the arc $\widehat{x_1x_2}$ we have that $\measuredangle x_2Ox_1=2\beta$. Repeating this process we obtain a sequence of points $x_3,$ $x_4, \ldots$ such that for every natural number $i$ it holds that $\measuredangle x_{i+1}Ox_i=2\beta$, and since $2\beta$ is an irrational multiple of $\pi$ we obtain a sequence of points $\{x_i\}_{i=1}^{\infty}$ which is dense in the circle $\mathbb S^1$. Moreover, since the distance from every support segment $[-x_i,x_{i+1}]$ to $O$ is equal to the distance from $[x_0,x_1]$ to $O$, we have that the support function of $K$ is constant and hence $K$ is a disc.
\item We suppose the theorem is proved for every number $n\leq r$, for some natural number $r$.
\item We will prove the theorem for $n=r+1$. Consider any point $x_0\in\mathbb S^r$ and let $\upsilon\in\mathbb S^r$ be any vector orthogonal to $x_0.$ Denote by $K_{\upsilon}:=K\cap \upsilon^{\bot}$ and $S_{\upsilon}:=\mathbb S^r\cap \upsilon^{\bot}$. Since $K_{\upsilon}$ have the hypothesis of the theorem from every point $x\in S_{\upsilon}$ we obtain that $K_{\upsilon}$ is a $r$-dimensional ball for every $\upsilon$ orthogonal to $x_0.$ We conclude that $K$ is a $(r+1)$-dimensional ball.
\end{enumerate}
Therefore, the proof is complete and $K$ is a ball. \fin

There are convex bodies in the plane, besides discs, which are contained in the interior of the unit circle $\mathbb S^1$ and such that all its projections onto $\mathbb S^1$ from points in $\mathbb S^1$ are congruent. These bodies indeed, are bodies which have the unit disc as one of is isoptic curves, i.e., curves under which the body is seen under a constant angle (see for instance \cite{Green}, \cite{Cies_Mier_Moz}). However, as we will see in Theorem \ref{proyecciones_congruentes_espacio}, if we consider projections over lines through the origin then the only convex body with the property of congruent projections is the disc. Moreover, the $(2n+1)$-dimensional version of this result, for every natural number $n$, still being true. Myroshnychenko had proved \cite{Myros} that two polytopes $K$ and $L$ are the same polytope if for any $x\in\mathbb S^{n-1}$ the projections of $K$ and $L$ onto $\mathbb S^{n-1}$ are congruent. 

We enunciate the following problem. 
 
\textbf{Problem 1.} It would be interesting to know whether or not a convex body $K\subset \mathbb R^n$ is a ball (or at least have centre of symmetry) if for every unit vector $x\in\mathbb S^{n-1}$ the sets $C(K,x)\cap x^{\bot}$ and $C(K,-x)\cap x^{\bot}$ are congruent. Even the case when both sets differ by a translation over the hyperplane $x^{\bot}$ seems to be interesting enough. The problem when we only know that $C(K,x)\cap x^{\bot}$ and $C(K,-x)\cap x^{\bot}$ have the same $(n-1)$-dimensional volume looks more difficult.

In the case of the Euclidean space with odd dimension we have the following.

\begin{theorem}\label{proyecciones_congruentes_espacio}
Let $K$ be a convex body contained in the interior of the unit ball of $\mathbb R^{2k+1}$, $k\geq 1$, and let $L$ be a $2k$-dimensional convex body. Suppose that for every unit vector $x\in\mathbb S^{2k}$ we have that $C(K,x)\cap x^{\bot}$ is a convex body congruent to $L$. Then $K$ is a Euclidean ball.
\end{theorem}

For the proof of this theorem we need to introduce an important notion: the concept of  \textit{field of congruent convex bodies}. Let $A\subset \mathbb{R}^n$ be an $n$-dimensional convex body, a \emph{field of bodies congruent to} $A$ is a continuous function $A(u)$ defined for $u$ in the unit sphere $\mathbb{S}^n=\{ x\in \mathbb{R}^{n+1}: \|x\|=1\}$, where $A(u)$ is a congruent copy of $A$ lying in a hyperplane of $\mathbb{R}^{n+1}$ perpendicular to $u$; here $A(u)$ is meant to be continuous in the Hausdorff metric. If in addition $A(u)= A(-u)$, for each $u$, we say $A(u)$ is a \textit{complete turning} of $A$ in $\mathbb{R}^{n+1}$. Clearly if all the $n$-dimensional sections of an $(n+1)$-dimensional convex body $K\subset \mathbb{R}^{n+1}$, through a fixed inner point are congruent, they give rise to a complete turning of some $n$-dimensional body $M\subset \mathbb{R}^{n}$ in $\mathbb{R}^{n+1}$.
 
In \cite{Mani} P. Mani proved  that \textit{if $n$ is even, fields of convex bodies tangent to $\mathbb{S}^n$} can be constructed only if $A$ is an $n$-dimensional ball. In \cite{Burton}, G. R. Burton used the methods of H. Hadwiger \cite{Hadwiger} and P. Mani \cite{Mani} to prove that a \textit{3-dimensional convex body $M\subset \mathbb R^4$ can be completely turned in $\mathbb R^4$ if and only if $M$ is centrally symmetric}. In \cite{Montejano1}, L. Montejano generalized Burton's result, i.e., \emph{convex bodies $M\subset \mathbb{R}^{n+1}$, $n \geq 2$, can be completely turned in $\mathbb{R}^{n+1}$ if and only if $M$ is centrally symmetric}. Montejano used his result to generalize a theorem of R. Schneider (see \cite{Schneider}) related to a famous problem of Banach.   
 
 
Let $I_{n+1}$ be the group of isometries of $\mathbb{R}^{n+1}$. Let us denote by $\Omega$ the space of all convex bodies in $\mathbb{R}^{n+1}$ with the topology induced by the Hausdorff metric. We consider the (continuous) map $\tau:\Omega \rightarrow \mathbb{R}^{n+1}$ which assigns to each $C\in \Omega$ the centre of its circumsphere. Let $C_1$ and $C_2\in \Omega$, we say that $C_1$ is \textit{congruent} to $C_2$ if there is $g\in I_{n+1}$ such that $g(C_1)=C_2$.

\emph{Proof of Theorem \ref{proyecciones_congruentes_espacio}.} We denote by $A(x)$ the set $C(K,x)\cap x^\perp$. The map $\alpha: \mathbb{S}^{2k} \rightarrow \Omega$  which assigns to each $x$ the set $A(x)$, is continuous. Indeed, if $x\in \mathbb{S}^{2k}$ and the sequence of points $\{x_r\}\subset \mathbb{S}^{2k}$ is such that $x_r \rightarrow x$, when $r \rightarrow \infty$, then 
$$S(K,x_r) \rightarrow S(K,x) , \textrm{ }  \textrm{ } \textrm{when} \textrm{ } \textrm{ } r \rightarrow \infty.$$
On the other side, $x^{\perp}_r \rightarrow x^{\perp}$, thus $A(x_r) \rightarrow A(x)$. 
Now we define the set $B(x)$ as $A(x)+(x-\tau(A(x)))$. Since $B(x)\subset x^{\perp}+x$, the family of convex bodies $\{B(x): x \in \mathbb{S}^{2k} \}$ varies continuously and since all the bodies $B(x)$ are congruent to $L$, we have found a field of bodies congruent to $L$. By Mani's theorem $L$ is a sphere. Now we apply a theorem due to G. Bianchi and P. Gruber (see \cite{Bianchi}) which affirm that if from any point $x$ in the boundary of a convex body $M$, for a convex body $N\subset \text{int}\, M$ it holds that $S(N,x)$ is an ellipsoidal cone, then $N$ is an ellipsoid. Now, is very easy to see that only Euclidean balls have the property stated in the theorem. Therefore, we conclude that $K$ is a Euclidean ball in $\mathbb R^{2k+1}$. \fin

The following theorem is the $2$-dimensional case of Theorem \ref{proyecciones_congruentes_espacio}, however, its proof cannot be obtained from the proof of Theorem \ref{proyecciones_congruentes_espacio}. For this reason, we include here in a separate statement.

\begin{theorem}\label{proyecciones_congruentes}
Let $K$ be a convex body contained in the interior of the unit disc of $\mathbb R^2$  such that for every unit vector $x\in\mathbb S^1$ we have that $C(K,x)\cap x^{\bot}$ is a segment of constant length. Then $K$ is a Euclidean disc.
\end{theorem}

\begin{figure}[H]
    \centering
    \includegraphics[width=0.6\textwidth]{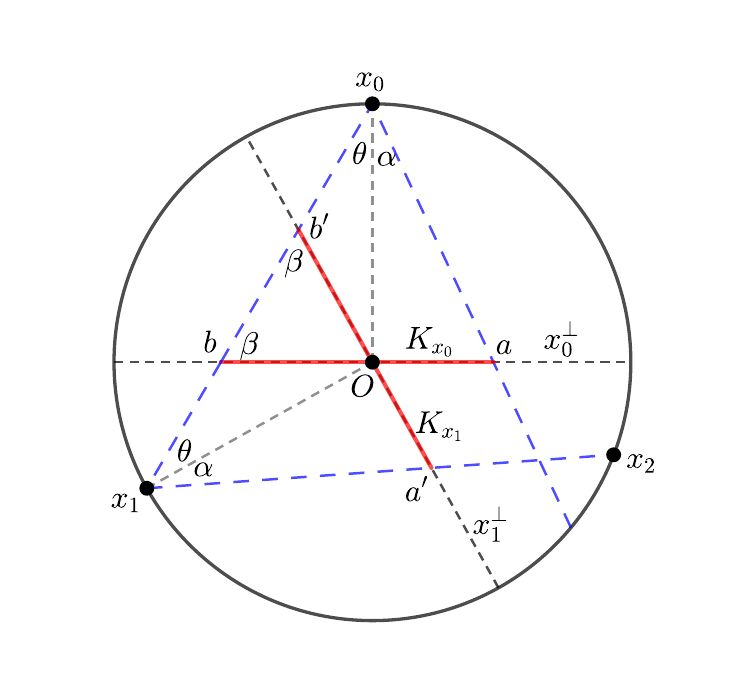}
    \caption{}
    \label{sucesion2}
\end{figure}

\emph{Proof.} Using a similar argument as the used in the proof of Theorem \ref{secciones_coinciden}, we obtain that the origin $O$ is contained in the interior of $K$. Suppose first that $K$ is seen under a constant angle $\varphi$ from every point $x_0\in\mathbb S^1$. Since the length of every segment $K_{x_0}=C(K,x_0)\cap x_0^{\bot}$ is also constant, we have that the two angles formed by the support lines of $K$ from $x_0$ with the segment $[x_0,O]$ are $\alpha$ and $\theta$, for two fixed values $\alpha$ and $\theta$ such that $\alpha+\theta=\varphi.$ Let $x_1\in\mathbb S^1$ be the point such that the segment $[x_1,x_0]$ is supporting $K$, as shown in Fig. \ref{sucesion2}. Let $b,a$ be points in the line $x_0^{\bot}$ such that $[b,a]=K_{x_0}$, and suppose with out loss of generality that $\measuredangle bx_0O=\theta$ and $\measuredangle ax_0O=\alpha$. Since the triangle $\triangle x_0Ox_1$ is isosceles, we have that $x_0^{\bot}$ and $x_1^{\bot}$ make equal angles with $[x_1,x_0]$. If $[b',a']$ is the projection of $K$ over the line $x_1^{\bot}$ from the point $x_1$ (as shown in Fig. \ref{sucesion2}), we have that $\measuredangle a'x_1O=\alpha$. Let $x_2$ be the point where the line $x_1a'$ intersects again to $\mathbb S^1.$ Continuing this process, we obtain a sequence of points, $x_0,$ $x_1$, $x_2$, $\ldots$, probably a finite sequence, such that the support function of $K$ takes only the values $\sin \alpha$ and $\sin \theta.$ No matter what is the initial point $x_0$, the associated sequence of points generates only the values $\sin\alpha$ and $\sin\theta$ for the support function. Since the support function of a convex body is a continuous function, we must have that $\theta=\alpha$ and hence $K$ is a disc.

Now, suppose the angle $\varphi$ is not the same for every $x\in\mathbb S^1$, then there is an interval $[\varphi_1,\varphi_2]$, such that $\varphi$ takes every value in this interval in a continuous way. Consider a point $x_0\in\mathbb S^1$ such that $K$ is seen under an angle $\varphi$ from $x_0$, with $\frac{\varphi}{\pi}$ an irrational number.  Again, suppose that $\measuredangle bx_0O=\theta$ and $\measuredangle ax_0O=\alpha$. This time the sequence of points $x_0$, $x_1,$ $x_2$, $\ldots$, is an infinite sequence, an moreover, the set of points $\{x_0,x_1,x_2, \ldots\}$ is dense in the circle $\mathbb S^1$. If follows that the support function of $K$ takes only the values $\sin\alpha$ and $\sin\theta$, and by the continuity of $p$, we must have that $\theta=\alpha.$ We conclude again, that $K$ is a disc centred at $O$. \fin

\begin{remark}
It would be interesting to know whether or not $K\subset\mathbb R^n$ is a ball if we only know that the $(n-1)$-dimensional volume of $C(K,x)\cap x^{\bot}$, for every $x\in\mathbb S^{n-1}$, is equal to a given constant number $\lambda$.
\end{remark}


Although the following result can be obtained by a theorem of G. Bianchi and P. Gruber, we provide a direct proof here which uses only geometric arguments. We believe this results goes towards finding of a proof for Conjecture \ref{penumbras_planas}.

\begin{theorem}\label{penumbras_planas}
Let $K,L\subset \mathbb R^n,$ $n\geq 3$, be two convex bodies with $L\subset \emph{int}\, K$. Suppose that for every $x\in\partial K$ the grazes $\Sigma(L,x)$ are $(n-1)$-dimensional ellipsoids. Then $L$ is an ellipsoid.
\end{theorem}

\emph{Proof.} A convex body $L\subset\mathbb R^n$, $n\geq 3$, is an ellipsoid if and only if there is a point $q\in\text{int}\, L$ and a fixed integer number $2\leq k\leq n-1$, such that all its $k$-dimensional sections through $q$ are ellipsoids. In particular, the above statement holds for $k=2$ and since the graze of $H\cap L$ from a point $x$, for every $3$-dimensional plane $H$, is precisely $H\cap \Sigma(L,x)$, we have that it is only need to prove the theorem for $n=3.$ For the point $q$ we may consider any point in the interior of $L$.

\begin{figure}[H]
    \centering
    \includegraphics[width=.96\textwidth]{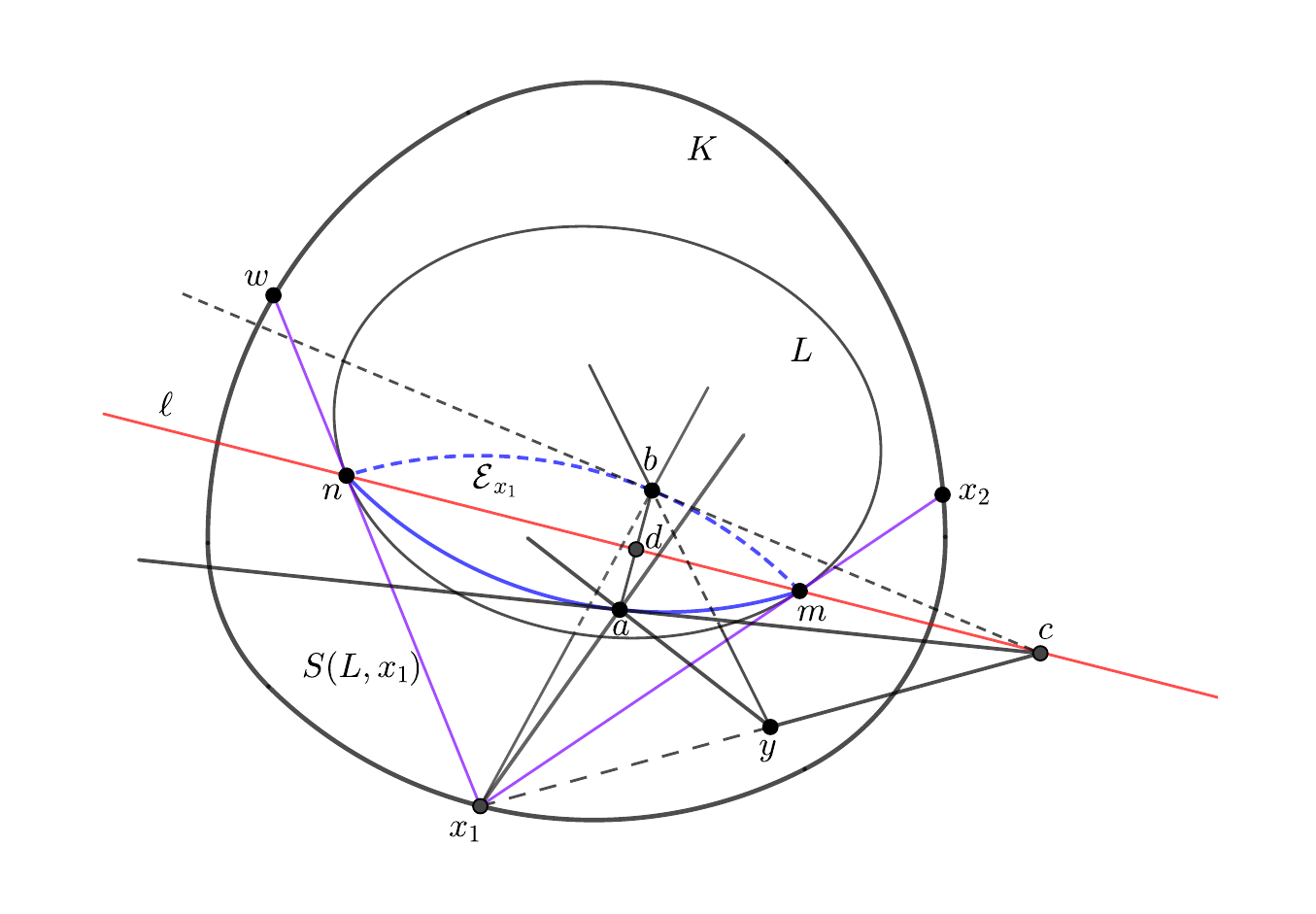}
    \caption{}
    \label{penumbra}
\end{figure}

In order to prove that $L$ is an ellipsoid we are going to prove that all the shadow boundaries of $L$ are planar, then we apply Blaschke's theorem (mentioned at the introduction) and the conclusion follows. Since ellipses are strictly convex it follows immediately that $L$ is strictly convex. For $x\in\partial K$ we denote by $\mathcal E_x$ the ellipse $\Sigma(L,x)$. Consider a unit vector $u\in\mathbb S^2$ and denote by $F_u$ the family of all grazes of $L$ contained in planes parallel to $u$. Let $x_1\in\partial K$ such that $\mathcal E_{x_1}\in F_u$ and let $\ell$ be the line of affine symmetry of $\mathcal E_{x_1}$ corresponding to $u$ (see Fig. \ref{penumbra}). Let $\Delta_1$ be the plane determined by $x_1$ and $\ell$. Consider any point $c\in \ell\setminus L$ and let $a$, $b$, the contact points of the supporting lines of $\mathcal E_{x_1}$ through $c$. Since $\ell$ is the line of affine symmetry of $\mathcal E_{x_1}$ corresponding to $u$, the segment $[a,b]$ is parallel to $u$ and its midpoint $d$ belongs to $\ell$. We denote by $y$ the point in $L(x,c)\cap\partial K$ ($L(x,c)$ denotes the line through $x$ and $y$), different of $x$, and denote by $\Omega$ the plane where $\mathcal E_y$ is contained. Since the planes $\text{aff}\{a,x_1,y\}$ and $\text{aff}\{b,x_1,y\}$ are supporting planes of $L$, which make contact with $L$ at $a$ and $b$, respectively, we have that $a$ and $b$ belong to $\mathcal E_y$. Thus $\Omega$ is parallel to $u$, i.e., $\mathcal E_y\in F_u$. Since $\Omega\cap\text{aff}\{a,x_1,y\}$ and $\Omega\cap\text{aff}\{b,x_1,y\}$ are supporting lines of $\mathcal E_y$ passing through a point $z\in L(x_1,y)$, the line of affine symmetry of $\mathcal E_y$ corresponding to $u$ is the line $L(z,d)$ which is contained in $\Delta_1$. 

Now, since $\mathcal E_{x_1}$ is an ellipse and the segment $[n,m]$ is the affine diameter of $\mathcal E_{x_1}$ contained in $\ell$, we have that the tangent lines to $\mathcal E_{x_1}$ at the points $n$ and $m$, are parallel to $u$. It follows that $n,m\in\Delta_1\cap\partial L$ are in the shadow boundary of $L$ with respect to the direction of illumination $u$, i.e., in $\Lambda(L,u)$. When we vary the point $c$ along the line $\ell$, we obtain that for every point $e$ of the arc $\widehat{x_2w}$ in $\Delta_1\cap \partial K$ (in clockwise sense, with respect to Fig. \ref{penumbra}), it holds that $\mathcal E_e\in F_u$ and consequently the corresponding line of affine symmetry of $\mathcal E_e$, with respect to $u$, intersects the boundary of $L$ in points that belong to $\Lambda(L,u)$. Again, we have that these points belong to $\Delta_1$. In this way we have that there is an arc in $\Delta_1\cap \partial L$ that belong to $\Lambda(L,u)$. 

\begin{figure}[H]
    \centering
    \includegraphics[width=.40\textwidth]{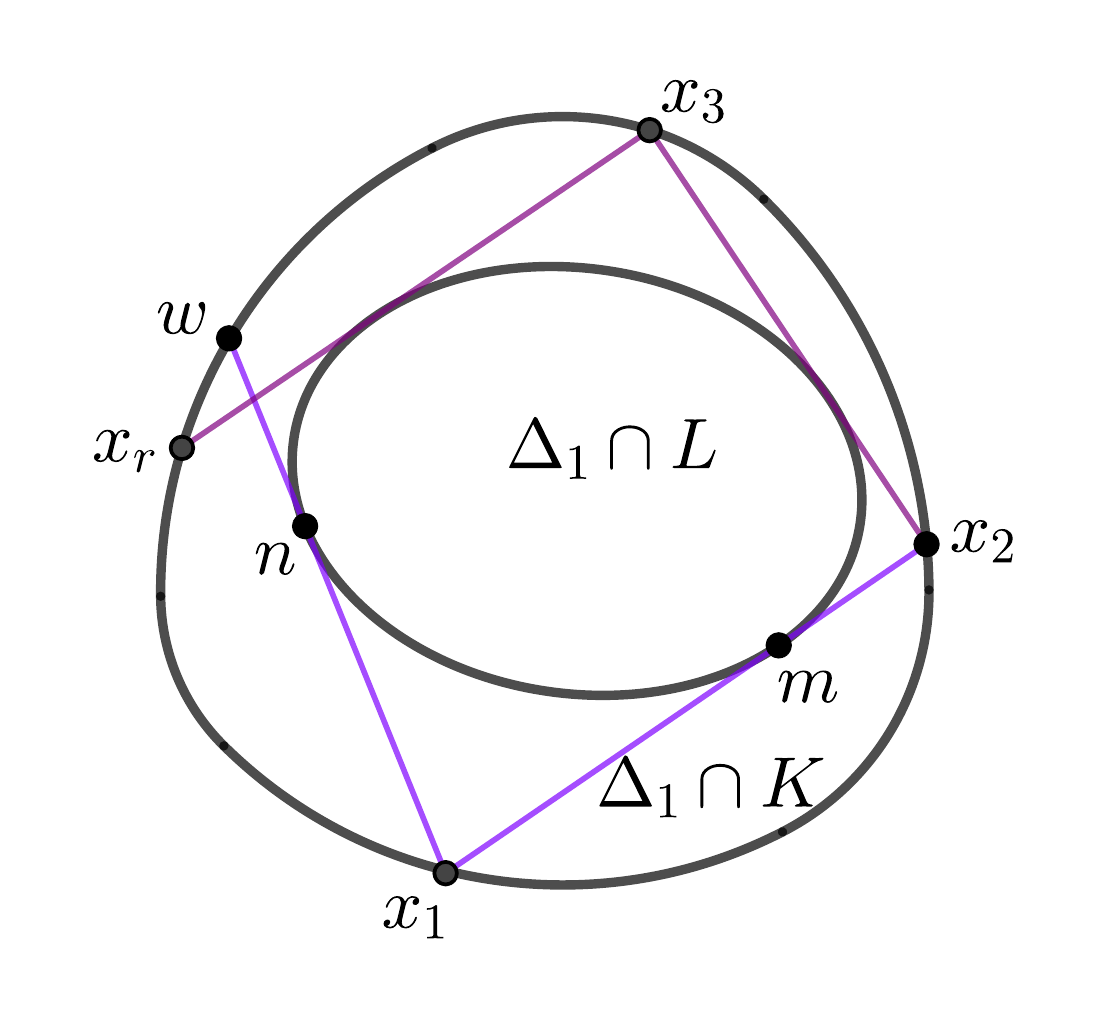}
    \caption{}
    \label{sucesion}
\end{figure}

To see that all the points in $\Delta_1\cap\partial L$ belong to $\Lambda(L,u)$, we need to prove that for every point $g$ in  the arc $\widehat{wx_2}$ (in clockwise sense) in $\Delta_1\cap\partial K$, it holds that $\mathcal E_g\in F_u$. In order to do this, we consider the sequence of points $x_2,$ $x_3$, $\ldots$ in $\Delta_1\cap\partial K$ (in counter-clockwise sense) such that every segment $[x_{i-1},x_i]$ is tangent to $\Delta_1\cap L$. Let $r$ be the minimum integer number such that $x_r$ is in the arc $\widehat{x_1w}$ in clockwise sense (see Fig. \ref{sucesion}). We have that for all the points $h$ in the arc $\widehat{x_2x_3}$, in counter-clockwise sense, it holds that $\mathcal E_h\in F_u$. We also have the same property for the points in the arc $\widehat{x_3x_4}$, and so on, until we arrive to the arc $\widehat{x_{r-1}x_r}$. The union of all these arcs cover $\Delta_1\cap\partial K$, hence, $\Delta_1\cap\partial L=\Lambda(L,u)$ (recall that $L$ is estrictly convex). This happens for every direction $u$, therefore, by Blaschke's theorem we conclude that $K$ is an ellipsoid. \fin

\begin{theorem}\label{demis}
Let $K,L \subset \Rn, n\geq 3$, be two convex bodies with $L \subset \emph{int} K $. Suppose that for every $x\in \partial K$ the sets $S(L,x) \cap L$ and $S(L,x) \cap K\setminus \{x\}$ are contained in parallel hyperplanes. Then $L$ and $K$ are homothetic and concentric ellipsoids.   
\end{theorem} 

\emph{Proof.} By the same argument given at the beginning of the proof of Theorem \ref{penumbras_planas}, we have that it is only need to prove the theorem for the case of dimension $n=3$. We will show that the body $L$ is strictly convex. Suppose, to the contrary, there exists a line segment $[a,b]\subset \partial L$. By virtue that $L\subset K$ there exist points $x,y\in \partial K$ such that $L(a,b)\cap \partial K=\{x,y\}$. By hypothesis there exist a plane $\Pi$ such that $S(L,x)\cap L\subset \Pi$. Since $[a,b]\subset S(L,x)\cap L$ it follows that $[a,b]\subset \Pi$. Consider a line $\tau\subset \Pi$ passing through $x$ such that $\tau\cap L=\emptyset$ and a point $q\in\text{int}\, L\setminus \Pi$. Now we rotate the plane 
$\text{aff}\{q,\tau\}$ with respect to the line $\tau$ until we get a plane $\Gamma$, 
$\Gamma \neq \Pi$, which is a supporting plane of $L$ and makes contact with $L$ at the point $s$. By the choice of $\tau$ and $\Gamma$, $s$ is not in $\Pi$. On the other hand, $s$ belongs to the set $S(L,x)\cap L$, that is, $s\in \Pi$. This contradiction shows that $L$ is strictly convex. 

Now, we are going to prove that for every supporting plane $\Pi$ of $K$ the section $\Pi \cap K$ is centrally symmetric with centre at $\Pi \cap L$. Then by Theorem O will follow that $L$ and $K$ are homothetic and concentric ellipsoids.

Let $\Pi$ be any supporting plane of $L$. In order to prove that $\Pi \cap K$ is centrally symmetric we will use the following well known characterization of central symmetry due to P. Hammer \cite{Hammer}:

\emph{Let $K\subset \mathbb R^n$, $n\geq 2$, be a convex body and $q\in \emph{int}\, K$. If for every line $l$ passing through $q$ there exist parallel supporting hyperplanes of $K$ at the points $l\cap \partial K$, then $K$ is centrally symmetric with centre at $q$.}

\begin{figure}[H]
    \centering
    \includegraphics[width=1.0\textwidth]{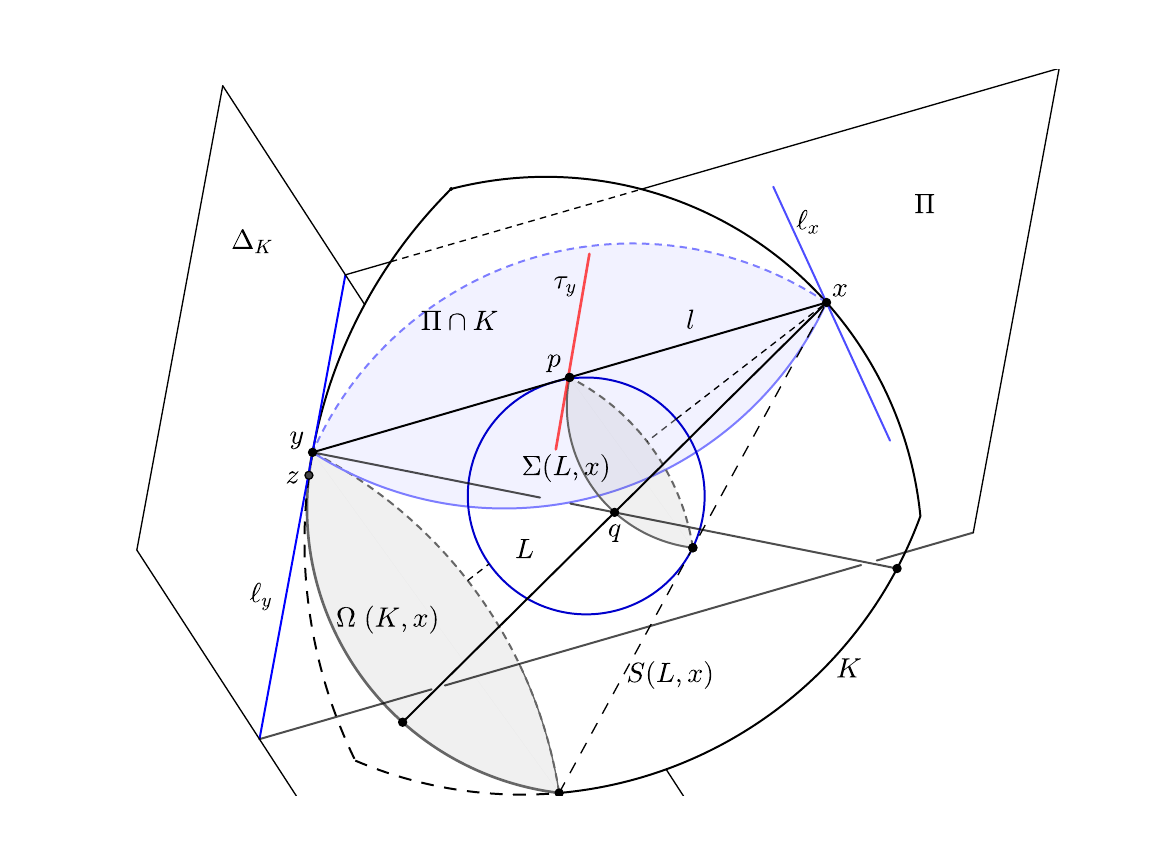}
    \caption{}
    \label{paralelas}
\end{figure}

Let $l\subset \Pi$ be a line passing through the point $p=\Pi \cap L$, we are going to prove that there are parallel supporting lines of $\Pi \cap K$ at the extreme points of the segment $l\cap \Pi \cap K$. We denote such points by $x,y$ and by $\Omega(K,x)$ the intersection $S(L,x)\cap K$ (See Fig. \ref{paralelas}). By hypothesis, there exists two parallel planes $\Delta_K$ and $\Delta_L$ such that $\Sigma(L,x) \subset \Delta_L$ and $\Omega(K,x) \subset \Delta_K$. Since $\Sigma(L,x)$ is a strictly convex curve and $\Omega(K,x)$ is homothetic to $\Sigma(L,x)$, with respect to the center of homothety $x$, we have that the line $\ell_y=\Pi \cap \Delta_K$ is a supporting line of $\Omega(K,x)$ at the point $y$ and intersects $\Omega(K,x)$ only at the point $y$. We have that $\ell_y$ must be also a supporting line of the section $\Pi\cap K$. Otherwise, let $z\in \Pi \cap K$ such that $z\in \ell_y$ and $z\neq y$, this imply that $z\in \ell_y\cap \Omega(K,x)$, which contradicts that $\ell_y$ touches $\Omega(K,x)$ only at $y$. Finally, we observe that $\ell_y$ is parallel to the line $\tau_y=\Pi \cap \Delta_L$. 

Let $\Delta_L'$ be the plane that contains to $\Sigma(L,y)$ and let $\Delta_K'$ be the plane parallel to $\Delta_L'$ through $x$. By the same argument as before, we have that there exists a supporting line $\ell_x$ of $\Pi \cap K$ at $x$ which is given by the intersection $\Pi \cap \Delta_K'$ and which is parallel to $\tau_x=\Pi \cap \Delta_L'$. If $\tau_x\neq \tau_y$ then there is another point of intersection between $\Sigma(L,x)$ and $\Sigma(L,y)$, $q\neq p$. Since $\Sigma(L,x)$ is a strictly convex curve, we have that $(p,q)\in\text{int}\, L$, however, since the lines $L(y,q)$ and $L(x,q)$ are supporting lines of $L$ through $q$, we have that the plane $\text{aff}\, \{y,q,x\}$ is a supporting plane of $L$ through $q$. This implies that $[p,q]\in\partial L$, which contradicts that $(p,q)\in\text{int}\, L$. Hence, $\tau_x=\tau_y$, i.e., $\ell_x$ is parallel to $\ell_y$. Since the line $l$ through $p$ was taken arbitrarily, we conclude by Hammer's theorem that $\Pi\cap K$ has centre of symmetry at $p$. Therefore, as we said before, $K$ and $L$ are homothetic and concentric ellipsoids. \fin

\end{document}